\documentclass{elsarticle}
\usepackage{mathtools,amsthm,amssymb}                           
\usepackage{physics} 
\usepackage[pagebackref]{hyperref}
\usepackage{colonequals}
\usepackage{theoremref}
\usepackage{enumitem}
\usepackage{fullpage}
\hypersetup{colorlinks = true,linkcolor = blue,anchorcolor =red,citecolor = blue,filecolor = red,urlcolor = red,pdfauthor=author}
\usepackage{bookmark}

\DeclareMathOperator{\conv}{\mathsf{conv}}
\DeclareMathOperator{\cl}{\mathsf{cl}}                  
\newcommand{\bracket}[1]{\left[#1 \right]}				
\DeclareMathOperator{\spec}{\mathsf{spec}}

\journal{Linear Algebra Appl.}

\newtheorem{theorem}{Theorem}[section]
\newtheorem{corollary}[theorem]{Corollary}
\newtheorem{lemma}[theorem]{Lemma}

\theoremstyle{definition}

\theoremstyle{remark}
\newtheorem{remark}[theorem]{Remark}

\numberwithin{equation}{section}

\begin{document}

\begin{frontmatter}
    \title{On a conjecture concerning totally extremal ideal Perron similarities}
    \author[add]{Erica J.~Artemis}
        \ead{ericaja@uw.edu}
    \author[add]{Pietro Paparella\corref{cor1}}
        \ead{pietrop@uw.edu}
        \affiliation[add]{organization={Division of Engineering \& Mathematics, University of Washington Bothell},
                    addressline={18115 Campus Way NE}, 
                    city={Bothell},
                    postcode={98011-8246}, 
                    state={WA},
                    country={U.S.A}}
        \cortext[cor1]{Corresponding author.}

\begin{abstract}
    Identifying ideal Perron similarities is a problem of central interest in the longstanding \emph{nonnegative inverse eigenvalue problem (NIEP)}. A normalized ideal Perron similarity is called \emph{totally extremal} if every entry has modulus one. Recently, Gershnik et al. [J.~Algebra {\bf 694} (2026), 782--800] proved that the character table of a finite Abelian group is totally extremal and conjectured the converse. 
    
    In this paper, we settle this conjecture in the affirmative by first showing that the rows of a totally extremal normalized ideal Perron similarity form a group under the Hadamard product. Then, it is shown that the rows of a nonsingular matrix form a group under the Hadamard product if and only if it is the character table of a finite Abelian group, and we further show that this group is isomorphic to the underlying group. These results extend the classical theorem due to Romanovsky and Karpelevi\v{c} on the unimodular eigenvalues of stochastic matrices in the complex unit disk to the setting of spectratopes in the unit ball of complex Euclidean space.
\end{abstract}


\end{frontmatter}

\section{Introduction}

Denote by $\mathsf{GL}_{n}$ the set of nonsingular complex matrices and by $e$ the all-ones vector of appropriate size. For $x \in \mathbb{C}^n$, denote by $D_x$ the diagonal matrix whose $k$-th diagonal entry is $x_k$. 

Given $S = [s_{ij}] \in \mathsf{GL}_n$, the set 
\[ 
\mathcal{C}(S) \colonequals \{ x \in \mathbb{C}^n \mid M_x = M_x(S) \colonequals SD_x S^{-1} \geqslant 0 \} \supseteq \{ \alpha e \mid \alpha \geqslant 0 \}, 
\] 
is a polyhedral cone~\cite[Theorem 4.6]{jp2025}, called the \emph{(Perron) spectracone of $S$}, and 
\[ 
\mathcal{P}(S) \colonequals \left\{ x \in \mathcal{C}(S) \mid M_x e = e \right\} \supseteq \{e\}, 
\]
is a polytope~\cite[Theorem 4.6]{jp2025}, called the \emph{(Perron) spectratope of $S$}. 

An entrywise nonnegative matrix $A$ is called \emph{stochastic} if every row sums to one, i.e., $Ae = e$. Since the spectral radius of a stochastic matrix is one, the polytope $\mathcal{P}(S)$ is contained in the unit ball $\mathbb{B}_\infty^n \colonequals \{ x \in \mathbb{C}^n \mid \| x\|_\infty \leqslant 1 \}$.

If there is a vector $x \in \mathbb{C}^n$ such that the matrix $M_x$ is irreducible (i.e., the \emph{directed graph} of $A$ is \emph{strongly connected} \cite[\S 3.2]{br1991}) and nonnegative, then $S$ is called a \emph{Perron similarity}. If $\mathcal{C}(S) = \mathcal{C}_r(S)$, where $\mathcal{C}_r(S)$ denotes the conical hull of the rows of $S$, then $S$ is called \emph{ideal}. It is possible for $S$ to be ideal but not a Perron similarity, and vice versa~\cite[Remarks 4.30 \& 4.31]{jp2025}. Identifying ideal Perron similarities remains a central focus in the study of the \emph{nonnegative inverse eigenvalue problem} (NIEP)~\cite{glp2026, jp2016, jp2017, jp2025}. 

If $S$ is an ideal Perron similarity, it can be assumed, without loss of generality, that $S e_1 = e$, 
\begin{equation}
    e_1^\top S = e^\top, \label{firstrowallones}    
\end{equation}
and $|s_{ij}| \leqslant 1$~\cite[pp.~8372--8375]{jp2025}. Such an ideal Perron similarity is called \emph{normalized}. As such, it is assumed throughout that every ideal Perron similarity is normalized. 

Recently, Gershnik et al.~\cite{glp2026} proved that the \emph{character table} of a finite group is an ideal Perron similarity. The case of a finite Abelian group is of particular interest because every entry of the corresponding matrix has unit modulus~\cite[p.~799]{glp2026}. Consequently, each row belongs to \( \partial\mathbb{B}_\infty^n = \mathbb{S}_\infty^n \colonequals \{ x \in \mathbb{C}^n \mid \|x\|_\infty = 1 \} \) and therefore lies on the boundary of the region $\mathbb{SL}^n$ comprising the spectra of stochastic matrices, viewed as vectors in $\mathbb{C}^n$~\cite[Section 3.2]{jp2025}. A normalized ideal Perron similarity is called \emph{totally extremal} if every entry has modulus one. Obviously, character tables of finite Abelian groups are totally extremal normalized ideal Perron similarities and characterizing the latter is of major interest for the NIEP.

To this end, Gershnik et al.~\cite{glp2026} conjectured that every totally extremal normalized ideal Perron similarity is the character table of a finite Abelian group. In this paper, we settle this conjecture in the affirmative. 

To see the significance of this result, the region in the complex plane 
\[
\Theta_n \colonequals \{ \lambda \in \mathbb{C} \mid \lambda \in \spec A,\ A \geqslant 0,\ Ae = e\},
\]
where $\spec A$ denotes the \emph{spectrum of $A$}, is contained in the unit disk $\mathbb{D} \colonequals \{ z \in \mathbb{C} \mid \vert{z}\vert \leqslant 1 \}$. Romanovsky~\cite{r1936}, and later Karpelevi\v{c}~\cite{k1951}, proved that 
\[
\Theta_n \cap \partial\mathbb{D} = \left\{ e^{\frac{2\pi}{q} \mathsf{i}} \mid p/q \in \mathcal{F}_n \right\},
\]
where $\mathcal{F}_n$ denotes the set of Farey fractions of order $n$. Consequently, our work extends Romanovsky's theorem for spectratopes to $\mathbb{B}_\infty^n$.

\section{Background \& Notation}

As usual, $\mathbb{N}$ denotes the set of positive integers, and for $n \in \mathbb{N}$, we let $\bracket{n} \colonequals \{ k \in \mathbb{N} \mid 1 \leqslant k \leqslant n\}$.

For a vector $x \in \mathbb{C}^n$, let $x_k$ (or $[x]_k$) denote its $k$-th entry. We say $x$ is \emph{totally nonzero} if $x_k \neq 0$ for all $k \in \bracket{n}$. The set of all such vectors is the \emph{algebraic torus} $(\mathbb{C}^\times)^n$. Under the Hadamard (entrywise) product $\circ$, the algebraic torus $(\mathbb{C}^\times)^n$ is an Abelian group with identity $e$.

Given $S \in \mathsf{GL}_n(\mathbb{C})$, $r_i \colonequals S^\top e_i \in \mathbb{C}^n$ denotes the $i$-th row of $S$ as a column vector and $R(S) \colonequals \{ r_i \mid i \in [n] \}$ denotes the set of row vectors of $S$. In particular, it follows from \eqref{firstrowallones} that $r_1 = e \in R(S)$ whenever $S$ is a normalized ideal Perron similarity. 

For a nonempty subset $X$ of $\mathbb{C}^n$, 
the \emph{convex hull of $X$}, denoted by $\conv{X} = \conv(X)$ is defined by  
\[ \conv {X} = \left\{ \sum_{k=1}^m \alpha_k x_k \in \mathbb{C}^n \mid m \in \mathbb{N},~x_k \in X,\ \sum_{k=1}^m \alpha_k = 1,\ \alpha_k \geqslant 0 \right\}. \]
The convex hull of a finite multiset $\{x_1,\ldots,x_n\}$ is abbreviated as $\conv(x_1,\ldots,x_n)$. Given $S \in \mathsf{GL}_n$, $\mathcal{P}_r(S) \colonequals R(S) = \conv (r_1, \ldots, r_n)$.

A subset $K$ of $\mathbb{C}^n$ is called \emph{convex} if $\conv(x,y) \subseteq K,\ \forall x,y \in K$. A point $x$ of a convex set $K$ is called an \emph{extreme point} (of $K$) if $u \notin K$ or $v \notin K$, whenever \(x = \alpha u + (1-\alpha) v\), with \(\alpha \in (0,1)\) and $u \ne v$, i.e., $x$ does not lie in any open line segment contained in $K$. 

Note for future reference that $\mathbb{B}_\infty^n$ is convex and any vector $v\in \mathbb{C}^n$ such that $\vert{v_k}\vert = 1,\ \forall k \in \bracket{n}$ is an extreme point.

\subsection{Characters \& Representations}

Henceforth, $G$ denotes a finite group with identity $e_G$. For ease of notation, the group operation $a \cdot b$ is abbreviated to $ab$. 

If $g \in G$, then $\cl g = \cl(g) \colonequals \{aga^{-1}\mid a \in G\}$ denotes the \emph{conjugacy class of $g$}. Because $\{ \cl g \mid g \in G \}$ forms a partition of $G$, we write $x \sim y$ whenever $x,y \in \cl g$.

A homomorphism $\rho:G \longrightarrow \mathsf{GL}_n(\mathbb{C})$ is called a \emph{(matrix) representation (of $G$)}. If $\rho$ is a representation, then $\rho(e_G) = I_n$ as $\rho$ is a homomorphism. 

If $\rho$ and $\sigma$ are representations, then:

\begin{itemize} 
    \item the \emph{Kronecker} or \emph{tensor product} of $\rho$ and $\sigma$, denoted by $\rho \otimes \sigma$, is the representation defined by $(\rho \otimes \sigma)(g) = \rho (g) \otimes \sigma (g),\ \forall g \in G$; and
    \item $\rho$ and $\sigma$ are called \emph{isomorphic} or \emph{similar}, denoted as $\rho \sim \sigma$, if there is an invertible matrix $S$ such that $\rho(g) = S\sigma(g) S^{-1},\ \forall g \in G$.
\end{itemize}

Given a representation $\rho$, the function $\chi_\rho : G \longrightarrow \mathbb{C}$ defined by 
\[ \chi_\rho(g) \colonequals \operatorname{tr}(\rho(g)),\ \forall g \in G, \] 
is called the \emph{character of $\rho$}. 
The \emph{dimension} (or \emph{degree}) of $\rho$, denoted by $\dim \rho$, is the quantity $\chi_\rho (e_G) = \operatorname{tr}(I_n) = n$. 
The dimension of $\chi_\rho$ refers to the same value. Because the trace is similarity invariant, equivalent representations possess identical characters; that is,
\begin{equation}
    \rho \sim \sigma \implies \chi_\rho(g) = \chi_\sigma(g)\ \forall g \in G. \label{simreps}
\end{equation}
In the special case where $\rho$ and $\sigma$ are one-dimensional, $\rho \sim \sigma$ if and only if $\rho(g) = \sigma(g)$ for all $g \in G$.

The following established properties (see, e.g., Fulton and Harris \cite[p.~13]{fh1991} or Serre \cite[Propositions 1 and 2]{s1977}) will be useful later:
\begin{align}
    \chi_{\rho \otimes \sigma} (g) &= \chi_\rho (g) \cdot \chi_\sigma (g)   \label{tensorprod} \\
    \chi_\rho (g) &= \chi_\rho (h),\ g \sim h.                               \label{classfun}
\end{align}

Given a representation $\rho$ of $G$ on $\mathbb{C}^n$, a subspace $W \subseteq \mathbb{C}^n$ is called \emph{$G$-invariant} if 
\[ \rho(g) W \subseteq W,\ \forall g \in G. \] 
The trivial subspaces $\{ 0 \}$ and $\mathbb{C}^n$ are always $G$-invariant. If $\mathbb{C}^n$ has a proper $G$-invariant subspace, then $\rho$ is called \emph{reducible} or a \emph{reducible representation}; otherwise, $\rho$ is called \emph{irreducible} or an \emph{irreducible representation}. Note that all one-dimensional representations are intrinsically irreducible. Furthermore, the tensor product of one-dimensional representations remains one-dimensional and is therefore also irreducible. A character corresponding to an irreducible representation is called an \emph{irreducible character}. The map $\rho_1 : G \longrightarrow \mathsf{GL}_n$, defined by $\rho_1(g) = I_n,\ \forall g \in G$, is called the ($n$-dimensional) \emph{trivial representation} and is irreducible if and only if $n = 1$.

The number of distinct irreducible representations (up to isomorphism) of $G$ equals the number of its distinct conjugacy classes (see, e.g., Artin \cite[Theorem 10.4.6(b)]{a1991}; Fulton and Harris \cite[Corollary 2.13]{fh1991}; or Serre \cite[Theorem 7]{s1977}). 

Suppose that $\rho_1, \ldots, \rho_n$ are the distinct irreducible representations (up to isomorphism) of $G$ and $\cl (g_1)$, $\ldots$, $\cl (g_n)$ are the distinct conjugacy classes of $G$. For ease of notation, we abbreviate $\chi_{\rho_k}$ to $\chi_k$. The $n$-by-$n$ matrix $Q = Q(G) = [q_{ij}] \colonequals [\chi_i(g_j)]$ is called the \emph{character table of $G$}. Properties \ref{simreps} and \ref{classfun} ensure that the matrix $Q$ is unique up to permutation of its rows or columns. 

The fundamental theorem of finite Abelian Groups asserts that if $G$ is a finite Abelian group, then there are prime numbers $p_1,\dots,p_d$ (not necessarily distinct) and positive integers $n_1,\ldots,n_d$ such that 
$$G \cong \bigoplus_{k=1}^d \mathbb{Z}_{p_k^{n_k}}.$$
Given that the \emph{discrete Fourier transform matrix} $F_m$ is the character table of $\mathbb{Z}_m$, it follows that 
$$Q(G) = \bigotimes_{k=1}^d F_{p_k^{n_k}},$$
where $\otimes$ denotes the matrix \emph{Kronecker} or \emph{tensor product}. Consequently, because discrete Fourier transform matrices are symmetric and the Kronecker product preserves symmetry, it can be assumed that the character table of any finite Abelian group $G$ is symmetric.

\section{Main Results}

While of independent interest, the following result will also be used to show that the rows of a totally extremal normalized ideal Perron similarity form a subgroup of the algebraic torus $(\mathbb{C}^\times)^n$. 

\begin{lemma}
    If $S$ is an ideal Perron similarity, then $r_i \circ r_j \in \mathcal{P}_r(S)$ for every $(i,j) \in \bracket{n}^2$. \thlabel{lem:hadprod}
\end{lemma}

\begin{proof}
    As $S$ is a normalized ideal Perron similarity, we have $\mathcal{P}(S) = \mathcal{P}_r(S)$ \cite[Theorem 4.22]{jp2025}. Thus, it suffices to show that $r_i \circ r_j \in \mathcal{P}(S),\ \forall (i,j) \in \bracket{n}^2$. 
    
    As $S$ is ideal, $r_i\circ r_j \in \mathcal{C}_r(S),\ \forall (i,j) \in \bracket{n}^2$ \cite[Theorem 3.5]{glp2026} and because $\mathcal{C}(S) = \mathcal{C}_r(S)$, we have $r_i\circ r_j \in \mathcal{C}(S),\ \forall (i,j) \in \bracket{n}^2$, i.e., $M_{r_i\circ r_j} \geqslant 0$. Since $Se_1 = e$, the first entry of $r_i\circ r_j$ is one because $[r_i]_1 = [r_j]_1 = 1$. Thus, 
    \begin{equation}
        M_{r_i\circ r_j} e = (S D_{r_i\circ r_j} S^{-1})(Se_1) 
        = \left(SD_{r_i\circ r_j} \right)e_1 
        = Se_1 
        = e,
    \end{equation}
    i.e., $M_{r_i\circ r_j}$ is stochastic and the result is established.
\end{proof}

\begin{theorem}
    If $S$ is a totally extremal ideal Perron similarity, then $\left(R(S),e,\circ\right) \leqslant (\mathbb{C}^\times)^n$. \thlabel{(i)implies(ii)}
\end{theorem}

\begin{proof}
    Recall that $e \in R(S)$; thus, by the finite subgroup test, it suffices to show that $R(S)$ is closed under the Hadamard product. To this end, because $S$ is an ideal Perron similarity, \thref{lem:hadprod} implies that $r_i \circ r_j \in \mathcal{P}_r(S)$ for every $(i,j) \in \bracket{n}^2$. Given that the entries of $r_{i}$ and $r_{j}$ have unit modulus, the same holds for the vector $r_i \circ r_j$. Consequently, $r_i \circ r_j$ is an extreme point of $\mathbb{B}_\infty^n$. Since the polytope $\mathcal{P}_r(S)$ is contained within $\mathbb{B}_\infty^n$, it follows that $r_i \circ r_j$ must also be an extreme point of $\mathcal{P}_r(S)$. Finally, because $r_i \circ r_j$ is a convex combination of the elements of $R(S)$---all of which are themselves extreme points---this can only occur if $r_i \circ r_j = r_k$ for some $k \in \bracket{n}$.
\end{proof}

\begin{theorem}
    If $Q \in \mathsf{GL}_n$, then $\left(R(Q),e,\circ\right) \leqslant (\mathbb{C}^\times)^n$ if and only if $Q$ is the character table of a finite Abelian group. \thlabel{(ii)iff(iii)}
\end{theorem}

\begin{proof}
For necessity, assume that $G \colonequals \left(R(Q),e,\circ\right) \leqslant (\mathbb{C}^\times)^n$. Note that $q_{ij} \ne 0$ since $r_k \in (\mathbb{C}^\times)^n$. Let $\hat{Q}$ be the character table of $G$, which can be assumed to be symmetric, and denote by $\hat{\rho}_1, \ldots, \hat{\rho}_n$ the distinct irreducible one-dimensional representations of $G$. 

For each $i,j \in \bracket{n}$, let $\rho_j \colon G \to\mathsf{GL}_1$ be the mapping defined by $\rho_j(r_i) = [q_{ij}]$. We demonstrate that $\rho_j$ is a representation: to this end, if $i, i' \in \bracket{n}$, then there is a positive integer $k \in \bracket{n}$ such that $r_i \circ r_{i'} = r_k$  since $G$ is a group. Thus,
\begin{equation}
    \rho_j(r_i \circ r_{i'}) = \rho_j(r_k) = [q_{kj}] = [q_{ij}q_{i'j}] = \rho_j(r_i)\rho_j(r_{i'}), \label{rhocalc}
\end{equation}
confirming that $\rho_j$ is a one-dimensional representation. Furthermore, these representations are distinct: indeed, if $\rho_j = \rho_{j'}$ for some $j \neq j'$, then 
\[ [q_{ij}] = \rho_j(r_i) = \rho_{j'}(r_i) = [q_{ij'}] \implies q_{ij} = q_{ij'},\ \forall i \in \bracket{n}, \]
which contradicts the invertibility of $Q$. Thus, there is a positive integer $j' \in \bracket{n}$ such that $\rho_j = \hat{\rho}_{j'}$. After relabeling, we may assume that $\rho_j = \hat{\rho}_{j},\ \forall j \in \bracket{n}$. Hence, 
$$q_{ij} = \tr \rho_j(r_i) = \tr \hat{\rho}_j(r_i) = \hat{q}_{ji} = \hat{q}_{ij},$$ 
i.e., $Q = \hat{Q}$.

For sufficiency, assume that $Q$ is the character table of a finite Abelian group $G$. By the finite subgroup test, it suffices to show non-emptiness and closure. 

Since $\rho_1$ is the irreducible trivial representation, note that $q_{1,j} = \rho_1(g_j) = 1,\ \forall j \in \bracket{n}$. Thus, $e_1^\top Q = e^\top$, i.e., $r_1 = e \in R(Q)$. Let $r_i$ and $r_j$ be rows of $Q$ corresponding to the characters $\chi_{\rho_i}$ and $\chi_{\rho_j}$, respectively. By Property~\ref{tensorprod}, the $k$-th entry of their Hadamard product satisfies
\[ [r_i \circ r_j]_k = \chi_{\rho_i}(g_k)\chi_{\rho_j}(g_k) = \chi_{\rho_i \otimes \rho_j}(g_k) \]
and because $\rho_i \otimes \rho_j$ is irreducible, its character must correspond to one of the rows of $Q$. Thus, $r_i \circ r_j \in R(S)$, and the result is established.
\end{proof}

\begin{corollary}
    For $Q \in \mathsf{GL}_n$, the following are equivalent: 
    \begin{enumerate}
    [label=\roman*)]
        \item \label{equivi} $Q$ is a totally extremal ideal Perron similarity;
        \item \label{equivii} $(R(Q),e,\circ) \leqslant (\mathbb{C}^\times)^n$; and
        \item \label{equiviii} $Q$ is the character table of a finite Abelian group.
    \end{enumerate} \thlabel{finitegroups}
\end{corollary}

\begin{proof}
    \thref{(i)implies(ii)} yields the implication \ref{equivi} $\implies$ \ref{equivii} and \thref{(ii)iff(iii)} establishes the biconditional \ref{equivii} $\iff$ \ref{equiviii}. Lastly, Gershnik et al.~\cite[\S 7]{glp2026} proved that \ref{equiviii} $\implies$ \ref{equivi}.
\end{proof}

\begin{remark}
    \thref{finitegroups} generalizes the classic well-known fact \cite[p.~314, Proposition 18]{df2004} that any finite subgroup of $\mathbb{C}^\times$ is of the form 
    $$\Omega_n \colonequals \left\{ e^{\frac{2\pi k}{n}\mathsf{i}} \mid 0 \leqslant k \leqslant n - 1 \right\}$$ for some $n \in \mathbb{N}$.
\end{remark}

Finally, we characterize the structure of the Abelian group $(R(Q), e, \circ)$.

\begin{theorem}
    If $Q\in \mathsf{GL}_n$ is the character table of a finite Abelian group $G$, then $(R(Q),e,\circ) \cong G$.
\end{theorem}

\begin{proof}
For definiteness, let $G = \{g_1, \dots, g_n\}$. The \emph{character group} of $G$, denoted by $\hat{G}$, is the group of its one-dimensional characters under pointwise multiplication and because $G$ is Abelian, we have $G \cong \hat{G} = \{ \chi_1, \dots, \chi_n\}$ \cite[Exercise 5.5, p~318]{a1991}. Thus, it suffices to show that $(R(Q), e, \circ) \cong \hat{G}$.

Since
    \[ r_i \circ r_{i'} = r_k \iff \chi_i(g_j) \chi_{i'}(g_j) = \chi_k(g_j),\ \forall j \in \bracket{n} \iff \chi_i\chi_{i'} = \chi_k, \]
it follows that the map $\varphi: R(Q) \to \hat{G}$ defined by 
$$\varphi(r_i) = \chi_i,\ \forall i \in \bracket{n}$$ 
is a homomorphism. The map $\varphi$ is clearly surjective and is injective because 
\[ \varphi(r_i) = \varphi(r_{i'}) \iff \chi_i = \chi_{i'} \iff \chi_i(g_j) = \chi_{i'}(g_j),\ \forall j \in \bracket{n} \iff r_i = r_{i'}. \qedhere \]
\end{proof}

\bibliographystyle{abbrv}
\bibliography{refs}

\begin{thebibliography}{10}

\bibitem{a1991}
M.~Artin.
\newblock {\em Algebra}.
\newblock Prentice Hall, Inc., Englewood Cliffs, NJ, 1991.

\bibitem{br1991}
R.~A. Brualdi and H.~J. Ryser.
\newblock {\em Combinatorial matrix theory}, volume~39 of {\em Encyclopedia of
  Mathematics and its Applications}.
\newblock Cambridge University Press, Cambridge, 1991.

\bibitem{df2004}
D.~S. Dummit and R.~M. Foote.
\newblock {\em Abstract algebra}.
\newblock John Wiley \& Sons, Inc., Hoboken, NJ, third edition, 2004.

\bibitem{fh1991}
W.~Fulton and J.~Harris.
\newblock {\em Representation theory}, volume 129 of {\em Graduate Texts in
  Mathematics}.
\newblock Springer-Verlag, New York, 1991.
\newblock A first course, Readings in Mathematics.

\bibitem{glp2026}
D.~Z. Gershnik, A.~J. Lewis, and P.~Paparella.
\newblock Character tables are ideal {P}erron similarities.
\newblock {\em J. Algebra}, 694:782--800, 2026.

\bibitem{jp2016}
C.~R. Johnson and P.~Paparella.
\newblock Perron spectratopes and the real nonnegative inverse eigenvalue
  problem.
\newblock {\em Linear Algebra Appl.}, 493:281--300, 2016.

\bibitem{jp2017}
C.~R. Johnson and P.~Paparella.
\newblock Row cones, {P}erron similarities, and nonnegative spectra.
\newblock {\em Linear Multilinear Algebra}, 65(10):2124--2130, 2017.

\bibitem{jp2025}
C.~R. Johnson and P.~Paparella.
\newblock Perron similarities and the nonnegative inverse eigenvalue problem.
\newblock {\em Trans. Amer. Math. Soc.}, 378(12):8361--8389, 2025.

\bibitem{k1951}
F.~I. Karpelevi\v{c}.
\newblock On the characteristic roots of matrices with nonnegative elements.
\newblock {\em Izvestiya Akad. Nauk SSSR. Ser. Mat.}, pages 361--383, 1951.

\bibitem{r1936}
V.~Romanovsky.
\newblock Recherches sur les cha\^{i}nes de {M}arkoff.
\newblock {\em Acta Math.}, 66(1):147--251, 1936.
\newblock Premier M\'{e}moire.

\bibitem{s1977}
J.-P. Serre.
\newblock {\em Linear representations of finite groups}, volume Vol. 42 of {\em
  Graduate Texts in Mathematics}.
\newblock Springer-Verlag, New York-Heidelberg, french edition, 1977.

\end{thebibliography}

\end{document}